\documentclass[a4paper,12pt]{article}

\usepackage[latin1]{inputenc}
\usepackage[english]{babel}

\usepackage{amssymb}
\usepackage{amsmath}
\usepackage{amsthm}
\usepackage{amscd}
\usepackage{mathtools}
\usepackage{enumitem}

\usepackage{textcomp}
\usepackage{layout}
\usepackage{color}

\renewcommand{\Pr}{\mathbb P} 
\newcommand{\Ex}{\mathbf E} 

\newcommand{\Poi}{\mathbf P} 
\newcommand{\Nd}{\mathcal N} 

\newcommand{\Poy}{\mathsf S} 
\newcommand{\Rp}{\mathsf R} 

\newcommand{\Pp}{\mathsf P} 
\newcommand{\GP}{\mathsf D} 

\newcommand{\B}{\mathcal{B}} 
\newcommand{\Bbd}{\mathcal{B}_{0}} 

\newcommand{\Gsig}{\mathcal G}

\newcommand{\MX}{\mathcal M(X)}

\newcommand{\Mpm}{\mathcal M^{\cdot\cdot}}
\newcommand{\MpmX}{\mathcal M^{\cdot\cdot}(X)}
\newcommand{\Mpf}{\mathcal M^{\cdot\cdot}_{f}}

\newcommand{\MpfN}{\mathcal M^{\cdot\cdot}_{f}(\N)} 
\newcommand{\MN}{\mathcal M(\N)}
\newcommand{\MfN}{\mathcal M_{\rm f}(\N)}
\newcommand{\MffN}{\mathcal M_{\rm ff}(\N)}

\newcommand{\N}{\mathbb N}

\newcommand{\R}{\mathbb R}

\renewcommand{\d}{\mathrm{d}}
\newcommand{\supp}{\operatorname{supp}}
\newcommand{\dom}{\operatorname{dom}}

\renewcommand{\exp}{\operatorname{exp}}
\newcommand{\e}{\operatorname{e}}

\newcommand{\ccdot}{\,\cdot\,}

\theoremstyle{definition}
\newtheorem{defdefinition}{Definition}[section]

\theoremstyle{plain}
\newtheorem{defsatz}[defdefinition]{Theorem}
\newtheorem{defsatzdef}[defdefinition]{Theorem and Definition}
\newtheorem{defprop}[defdefinition]{Proposition}
\newtheorem{deflemma}[defdefinition]{Lemma}
\newtheorem{deffolgerung}[defdefinition]{Corollary}
\newtheorem{defbeispiel}[defdefinition]{Example}

\theoremstyle{remark}
\newtheorem{defbemerkung}[defdefinition]{Remark}

\newcommand{\prop}[1]{\begin{defprop}#1\end{defprop}}

\newcommand{\bem}[1]{\begin{defbemerkung}#1\end{defbemerkung}}
\newcommand{\lemma}[1]{\begin{deflemma}#1\end{deflemma}}
\newcommand{\lemman}[2]{\begin{deflemma}[#1]#2\end{deflemma}}
\newcommand{\korollar}[1]{\begin{deffolgerung}#1\end{deffolgerung}}

\newcommand{\C}{{\cal C}}

\newcommand{\A}{{\cal A}}

\title{A hydrodynamic limit and fluctuations for a Chinese restaurant-like process}

\author{Mathias Rafler\footnote{email {\sf rafler@ma.tum.de}}}

\begin{document}

\maketitle

\begin{abstract}
We study a Markov process constructed from the P\'olya sum process, which yields a kind of spatial version of the Chinese restaurant process, where each 'table' is assigned a 'location'. This construction firstly allows a definition of summation of independent processes, and secondly a monotone coupling for different parameters. Moreover, the process is related to a multi-particle random walk on the positive integers and its hydrodynamic limit and fluctuations are discussed.

\end{abstract}

\section{Introduction}

The chinese restaurant process, widely discussed e.g. in~\cite{PY97}, is a Markov process constructing a growing sequence of partitions identifying each block with the subset of customers at some restaurant sitting at one table. Starting with some number $n$ of customers in a restaurant sitting at some $k\leq n$ tables, a new customer chooses to sit at an empty table with probability proportional to 1 and at an already occupied table with probability propotional to the number of guests sitting at that table. That way the discrimination of customers according to their table forms a partition on the set of customers present at the restaurant. In~\cite{EGW12} inter alia the Martin boundary of the Chinese restaurant process is computed.

A closely related point process is the P\'olya sum process, the Papangelou process whose Papangelou kernel is given by
\begin{equation*}
  \pi(\mu,\ccdot)=z(\rho+\mu)
\end{equation*}
for some $z\in(0,1)$ and $\rho\in\MX$. Here new points are introduced given a configuration of points $\mu$ at a rate proportional to $\rho+\mu$, that is, a new point may either choose a location according to a fixed measure $\rho$, or join the location of a present point.

The aim is to construct a Markov process whose one-dimensional marginals are P\'olya sum processes, and whose paths are increasing. In contrast to the chinese restaurant process, points are not introduced at discrete times but continuously during the time intervall $[0,1)$. New `tables' are opened according to the measure $\rho$ which may be an infinte but locally measure. Thus during a time intervall $(t,t+h]$ an infinite but locally finite number of customers may arrive and choose to join an existing table or open a new one. In~\cite{mR13coxexit} this process was constructed and its exit space was characterized and thereby this process identified as a Gibbs process. We recall the definition and the basic properties of the P\'olya sum process section~\ref{sect:polya} and give the definition and properties of the spatial version of the Chinese restaurant-like process as well as its relation to the classical Chinese restaurant process (CRP) in section~\ref{sect:scrp}. 

An immediate consequence of using the spatial process is possibility of the definition of the sum of two independent CRP, which is again a CRP where the parameter is just the sum of the two parameters. Moreover, one obtains a coupling of CRP for different parameters in a monotone way such that at each time, one of the processes contains the partitions of the other one.

Section~\ref{sect:mprw} is dedicated to the relation of the Markov process to a random walk of particles on the positive integers: Considering each table as a particle nd the number of customers as its location, each particle hops to the right in the moment a new customer arrives at the corresponding table while new particles are introduced at 1 when a customer opens a new table. The dynamics is such that the numbers of particles at different sites are still independent. The basic tool will be state space transformations for Markov processes.

The question that arises is how this random walk behaves for large parameter, or how does the spatial chinese restaurant process behave in large domains. We study first the hydrodynmic limit in section~\ref{sect:hydra}, that is the non-random limiting process and the corresponding dynamics. Secondly, we identify the fluctuations around this mean and study their dynamics. The surprising result is that despite that they are independent for different sites, they do not evolve independently.

\section{The P\'olya sum process \label{sect:polya}}

Throughout the paper let $X$ be a Polish space with Borel $\sigma$-field $\B$ and a ring of bounded Borel sets $\Bbd$. By $F$ we denote the set of continuous functions from $X$ to $\R$, by $F_+\subset F$, $F_b\subset F$ the sets of non-negative, continuous functions and the set of continuous functions with bounded support, respectively. $F_{+,b}=F_+\cap F_b$.

When equipped with the vague topology, the space of locally finite measures $\MX$ is Polish as well as its closed subset of locally finite point measures $\MpmX$. Its Borel $\sigma$-field $\Gsig$ is generated by the evaluation mappings $\zeta_B:\MX\to\R_+$, $\zeta_B\mu=\mu(B)$. Let $\Gsig_B=\sigma(\zeta_{B'}:B\in\B, B\subset B)$. Often $\int f\d\mu$ is denoted by $\mu(f)$ or $\zeta_f\mu$. The elements of $\MX$ are partially ordered via $\nu\leq\mu$ iff this inequality holds for all non-negative test functions on $X$. For point processes this means that the point masses of $\nu$ are dominated by the ones of $\mu$ or, equivalently, that $\mu-\nu$ is still a point measure.

For measures on $\N$, denote by $\MfN$ the set of finite measures and by $\MffN$ the set on measures on $\N$ with finite first moment. Moreover, we need the notion of Gateaux differentiability of a function on these measure spaces: A function $f$ is Gateaux differentiable at $\nu$ in direction $\kappa$, if the limit 
\begin{equation*}
	\lim_{h\to0}\frac{1}{h}\bigl[f(\nu+h\kappa)-f(\nu)\bigr]
\end{equation*}
exists; in this case the limit is denoted by $f'(\nu)[\kappa]$. If this holds for all $\nu$ and $\kappa$, $f$ is said to be Gateaux differentiable. By iteration one obtains higher derivatives, and moreover a Taylor expnsion is available. See e.g.~\cite{rH82}.

A random measure is a probability measure on $\MX$, a point process is a random measure which is concentrated on $\MpmX$. For a random measure $\Rp$ the Campbell measure is defined as
\begin{equation*}
	C_\Rp(h)=\iint h(x,\nu)\nu(\d x)\Rp(\d\nu),
\end{equation*}
where $h:X\times\MX\to\R$ is measurable and non-negative or integrable. Point processes $\Pp$ satisfying some absolute continuity condition admit a particular disintegration of their Campbell measure,
\begin{equation*}
	C_\Pp(h)=\iint h(x,\mu+\delta_x)\pi(\mu,\d x)\Pp(\d\nu),
\end{equation*}
see e.g.~\cite{MWM79,NZ79,oK83,MKM78}. In this case, $\pi$ is called Papangelou kernel and $\Pp$ Papangelou process. Well known is the Poisson process $\Poi_\rho$ given by the kernel $\pi(\mu,\ccdot)=\rho$, $\rho\in\MX$. Of further special interest will be the P\'olya sum process $\Poy_{z,\rho}$ and the P\'olya difference process $\GP_{z,\rho}$ which are characterized by the sum kernel $\pi(\mu,\ccdot)=z(\rho+\mu)$, $z\in(0,1)$, $\rho\in\MX$, and the difference kernel $\pi(\mu,\ccdot)=z(\rho-\mu)$, $z>0$, $\rho\in\MpmX$, $\mu\leq\rho$. The latter processes were introduced and constructed in~\cite{hZ09,NZ11}.

We need two main important properties of the P\'olya sum process: Firstly the independence of increments, and secondly the convolution property $\Poy_{z,\rho+\nu}=\Poy_{z,\rho}\ast\Poy_{z,\nu}$.

The first aim is to construct a Markov process with increasing paths such that at time $t$ its one-dimensional distribution is $\Poy_{t,\rho}$. Main tools are two thinning and splitting given in~\cite{mR13coxexit}, which are adapted from~\cite{bN12}. 

\lemman{Sampling from P\'olya sum processes}{ \label{thm:thinnings:polyathinning}
Let $\Pp=\Gamma_q(\Poy_{z,\rho})=\int\GP_{\frac{q}{1-q},\mu}\Poy_{z,\rho}(\d\mu)$. Then $\Pp=\Poy_{\gamma,\rho}$, where $\gamma=\gamma(z,q)=\tfrac{zq}{1-z(1-q)}$.
}%
Via splitting a P\'olya sum process with higher density from one with lower density such that the former dominates the latter is constructed.

\lemman{Condensations of P\'olya sum processes}{ \label{thm:splittings:polyasplitting}
Let $0\leq\gamma\leq z<1$ and 
\begin{equation*}
  \Pp(\phi)=\iint \phi(\mu+\nu) \Poy_{\frac{z-\gamma}{1-\gamma},\rho+\nu}(\d\mu) \Poy_{\gamma,\rho}(\d\nu).
\end{equation*}
Then $\Pp=\Poy_{z,\rho}$.
}%

\section{The spatial Chinese Restaurant-like process \label{sect:scrp}}

Define $\Omega=\{\omega:[0,1)\to\MpmX: \omega\text{ cadlag, }\omega_s\leq\omega_t, s\leq t\}$ and equip $\Omega$ with the Skorohod topology. Define $Y_t(\omega)=\omega_t$. Since the state space $\MpmX$ is separable and complete, the following equations define a Markov process uniquely on $\Omega$,~\cite[Thm 4.1.1]{EK05}.
\begin{align*}
	Y_0&=0\\
	\Pr(Y_t-Y_s\in\ccdot|Y_s) &= \Poy_{\frac{t-s}{1-s},\rho+Y_s},\qquad 0\leq s<t<1.
\end{align*}
The transition kernels then are 
\begin{equation*}
	p_{s,t}(\nu,\ccdot)=\Poy_{\frac{t-s}{1-s},\rho+\nu}\ast\Delta_\nu
\end{equation*}
where $\Delta_\nu$ is the point process which realizes the configuration $\nu$ a.s. Because of Lemma~\ref{thm:splittings:polyasplitting}, $Y_t$ is $\Poy_{t,\rho}$-distributed for each $t$. The backward dynamics $p^\ast$ given by independent thinning,
\begin{equation*}
	p^\ast_{s,t}(\nu,\ccdot)=\Gamma_{\frac{s(1-t)}{t(1-s)}}(\nu)
\end{equation*}
for $s<t$ (one starts at time $t$ with configuration $\nu$ and goes back to times $s$). By the sampling lemma this means that 
\begin{equation*}
	p^\ast_{s,t}(\nu,\phi)=\GP_{\frac{s(1-t)}{t-s},\nu}(\phi).
\end{equation*}

\bem{ \label{bem:scrp:singleregion}
For fixed $B\in\Bbd$, $\bigl(Y_t(B)\bigr)_t$ is a Markov process with values in $\N_0$ such that for each $t$, $Y_t(B)\sim NB\bigl(\rho(B),t\bigr)$. The transition from $Y_s(B)$ to $Y_t(B)$ means to add a $NB\bigl(\rho(B)+Y_s(B),\tfrac{t-s}{1-s}\bigr)$ distributed random variable, and backwards to take $Y_s(B)$ conditioned on $Y_t(B)$ to be $B\bigl(Y_t(B),\frac{s(1-t)}{t(1-s)}\bigr)$ distributed.
}

\subsection{Relation to the Chinese Restaurant Process and arrival times}

We show in which way to obtain the chinese restaurant process from the spatial construction. This means firstly to go from continuous to discrete time and secondly to forget about spatial relations. Observations are only made in a bounded observation window.

For a bounded set $B$ and positive integer $m$ let
\begin{equation*}
	\tau_{B,m}=\inf\{t\geq 0:Y_t(B)\geq m\}
\end{equation*}
the time of the arrival of the $m$-th point in $B$. 

The following properties follow immediatly: If $\rho(B)>0$, then by remark~\ref{bem:scrp:singleregion}, $\tau_{B,m}<1$ a.s. for each $m$ since $Y_t(B)\to\infty$ a.s, and moreover, as $m\to\infty$, $\tau_{B,m}\to 1$ since $Y_t(B)$ is a.s. finite for $t<1$. We will see later that at each time instant only one point may arrive, i.e. a.s. we have $0<\tau_{B,1}<\tau_{B,2}<\ldots$.

Note that we get explicitly the distributions of the arrival times,
\lemma{ \label{thm:crp:arrivaltimedistr}
Let $B\in\B$ such that $0<\rho(B)<\infty$. Then for each non-negative integer $m$, $\tau_{B,m+1}\sim\beta\bigl(m+1,\rho(B)\bigr)$.
}
\begin{proof}
\begin{align*}
	\Pr(\tau_{B,m+1}> t) &= p_{0,t}(0,\zeta_B\leq m)=\Poy_{t,\rho}(\zeta_B\leq m)=(1-t)^{\rho(B)}\sum_{k=0}^{m}\frac{t^k}{k!}\rho(B)^{[k]},
\end{align*}
and differentiating yields the density of $\tau_{B,m}$
\begin{equation*}
	\frac{\rho(B)^{[m+1]}}{m!}t^m(1-t)^{\rho(B)-1}.\qedhere
\end{equation*}
\end{proof}

\bem{
There is the freedom to proof the lemma via the Cox representation of the P\'olya sum process. In this case, the dependence of $t$ is only due to the directing measure, a Poisson-Gamma process, see e.g.~\cite{mR13coxexit}: $\zeta_B~\sim\Gamma\bigl(t,\rho(B)\bigr)$.
}
The conditional distribution of the waiting times is again stronly connected to the beta distribution.
\lemma{
Conditioned on $\tau_{B,m}=s\in(0,1)$ with $B$ chosen as in Lemma~\ref{thm:crp:arrivaltimedistr} and $m\in\N$,
\begin{equation*}
	\frac{\tau_{B,m+1}-s}{1-s}~\sim\beta\bigl(1,\rho(B)+m\bigr).
\end{equation*}
}
Thus the waiting times can be obtained via the following procedure: Let $U_k\sim \beta\bigl(1,\rho(B)+k\bigr)$, $k\in\N$ be independent, then the first waiting time is $U_1$ and the $m$-th waiting time is $(1-U_1)\cdots(1-U_{m-1})U_m$. 
\begin{proof}
On the set $\{\tau_{B,m}=s\}\subset\{Y_s(B)=m\}$ we have
\begin{align*}
	\Pr(\tau_{B,m+1}>t|\tau_{B,m}=s) &= \int 1_{\{\zeta_B\leq m\}}(\mu)p_{s,t}(Y_s,\d\mu)\\
		&= \int 1_{\{\zeta_B=0\}}(\mu-Y_s)p_{s,t}(Y_s,\d\mu)\\
		&= \Poy_{\frac{t-s}{1-s},\rho+Y_s}(\zeta_B=0) = \left(1-\frac{t-s}{1-s}\right)^{\rho(B)+m}
\end{align*}
for $0<s<t<1$. Thus under the given scaling we observe the claimed $\beta$-distribution. \end{proof}

Connecting these two lemmas, we get the joint density function of the two successive arrival times $\tau_{B,m}$ and $\tau_{B,m+1}$ as
\begin{equation*}
	\frac{\rho(B)^{[m+1]}}{(m-1)!} \cdot \frac{(1-t)^{\rho(B)+m-1} s^{m-1}}{(1-s)^{m+1}} 1_{0\leq s\leq t\leq 1}.
\end{equation*}

\bem{
Assume that $\rho$ is a diffuse measure in the following. 
\begin{enumerate}
	\item Considering the Markov chain $(Y_{\tau_{B,m}})_m$, we construct random partitions as follows: At time $m=1$ the first point arrives at some $x_1\in B$ and restricted to $B$, $Y_{\tau_{B,1}}=y_{1,1}\delta_{x_1}$. At time $m-1$, again restricted to $B$, $Y_{\tau_{B,m-1}}=y_{m-1,1}\delta_{x_1}+\ldots+y_{m-1,k}\delta_{x_k}$ for some $k\leq m$, pairwise distinct $x_1,\ldots,x_k$ and positive integers $y_1,\ldots,y_k$ which sum up to $m-1$. The $m$-th point may arrive either at some point $x_{k+1}$ distinct from all present points with probability proportional to $\rho(B)$ and giving rise to $y_{m,k+1}=1$ and $y_{m,j}=y_{m-1,j}$, $j\leq k$; or it may arrive at some $x_j$, $j\leq k$ yielding $y_{m,j}=y_{m-1,j}+1$ and keeping all the other points. This is the construction of the chinese restaurant process with weight $\rho(B)$.
	\item The construction allows algebraic manipulations of chinese restaurant processes: Joining two independent processes yields again a chinese restaurant process with the new parameter being the sum of the old ones. More precisely, the two independent sequences of arrival times are ordered, and then the partitions are built up in the just descibed way. 
	
	In continuous time one just joins the partitions at each time. This even holds for infinite sums as long as the parameters are summable (i.e. one stays in a bounded set).
	\item By observing the process in two bounded sets $B'\subset B$, one obtains a monotone coupling in the sense that at each time, the partition of the process corresponding to the smaller set is contained in the partition of the process of the bigger set.
\end{enumerate}
}

\subsection{The generator}

We compute the generator of the (non-homogeneous) Markov process in the case of finite point processes,
\begin{equation*}
	\A_s\phi=\lim_{h\to 0} \frac{1}{h}\left[p_{s,s+h}(\ccdot,\phi)-\phi\right].
\end{equation*}
\lemma{ \label{thm:gen:gen}
Let $\rho$ and $\nu$ be locally finite measures and $\phi:\MX\to\R$ be bounded and $\Gsig_B$-measurable. Then
\begin{equation*}
	\A_s\phi(\nu) = \frac{1}{1-s} \int \phi(\nu+\delta_x)-\phi(\nu)\bigl(\rho+\nu\bigr)(\d x)
\end{equation*}
for $s<1$.
}
\begin{proof}
Sorting the summands, we obtain since $\phi$ is $\Gsig_B$-measurable,
\begin{align*}
	p_{s,s+h}(\nu,\phi)&-\phi(\nu) \\
		&= \begin{multlined}[t]
				\left[\Poy_{\frac{h}{1-s},\rho+\nu}(\zeta_B=0)-1\right]\phi(\nu)\\
				+ \int_{\zeta_B=1} \phi(\mu+\nu) \Poy_{\frac{h}{1-s},\rho+\nu}(\d\mu) 
				+\int_{\zeta_B\geq 2} \phi(\mu+\nu) \Poy_{\frac{h}{1-s},\rho+\nu}(\d\mu)
			\end{multlined}\\
		&= \begin{multlined}[t]
				\left[\left(1-\frac{h}{1-s}\right)^{\rho(B)+\nu(B)}-1\right]\phi(\nu)\\
				+ \left(1-\frac{h}{1-s}\right)^{\rho(B)+\nu(B)}\cdot \frac{h}{1-s}\int \phi(\mu+\delta_x) \bigl(\rho+\nu\bigr)(\d x)\\
				+ \int_{\zeta_B\geq 2} \phi(\mu+\nu) \Poy_{\frac{h}{1-s},\rho+\nu}(\d\mu).
			\end{multlined}
\end{align*}
Dividing the equation by $h$, the first summand tends to $-\tfrac{\rho(B)+\nu(B)}{1-s}\phi(\nu)$ as $h\to 0$. Since by the boundedness of $\phi$ the last one is of order $O(h)$, we get 
\begin{equation*}
	\A_s\phi(\nu) = \frac{1}{1-s}\left[ \int \phi(\nu+\delta_x)\bigl(\rho+\nu\bigr)(\d x) -\bigl(\rho(B)+\nu(B)\bigr)\phi(\nu) \right]
\end{equation*}
and do last arrangements.
\end{proof}

\bem{
One confirms by a direct computation Kolmogorov's equation, i.e. for some $\Gsig_B$-measurable, bounded  $\phi$,
\begin{equation*}
	\frac{\d}{\d t}\Poy_{t,\rho}(\phi)= \Poy_{t,\rho}(\A_t\phi).
\end{equation*}
}
\bem{
One may introduce a time-changed version $X_u=Y_{\frac{u}{1+u}}$, i.e. $Y_s=X_{\frac{s}{1-s}}$. Since the law of $X_u$ then is $\Poy_{\frac{u}{1+u},\rho}$, we get that the intensity measure $\Ex X_u$ of $X_u$ is $u\rho$, and moreover
\begin{equation*}
	\Pr(X_u\in\ccdot|X_v)=\Poy_{\frac{u-v}{1+u},\rho+X_v}
\end{equation*}
for $0\leq v<u$. Anyways, this time change does not make the family of generators time homogeneous.
}

\section{A multi-particle random walk representation \label{sect:mprw}}

The reduction of the spatial information admits a further point of view in terms of particles moving on the natural numbers always to the right while being introduced at 1. The basis is Dynkin's state space transformation theorem for Markov processes~\cite{eD65}: Roughly, if in this situation $\gamma:\MpmX\to\MpfN$ is a well-behaved function such that the transition kernels satisfy
\begin{equation*}
	p_{s,t}(\mu_1,f\circ\gamma)=p_{s,t}(\mu_2,f\circ\gamma)
\end{equation*}
for all $0\leq s<t<1$, and $\mu_1,\mu_2\in\MpmX$ such that $\gamma\mu_1=\gamma\mu_2$, and if one defines $\gamma^\ast f=f\circ\gamma$, then $\gamma(Y_t)$ is a Markov process with transition kernel given by
\begin{equation*}
	p^\gamma_{s,t}(\eta,f)=p_{s,t}(\mu,\gamma^\ast f)
\end{equation*}
for some suitable $\mu\in\{\gamma=\eta\}$, the transition semigroup is given by $\gamma^\ast T^\gamma_{s,t}=T_{s,t}\gamma^\ast$ and its infinitesimal generator by $\gamma^\ast \A^\gamma_{s,t}=\A_{s,t}\gamma^\ast$.

Again we assume $\rho$ to be a diffuse measure. Fix some bounded set $B$ and denote by the random variable $U_{B}(j)$ the number of points in $B$ of multiplicity $j$ of some point configuration in $\MpmX$. By $U_{B}\in\Mpf(\N)$ denote the measure with weights $U_{B}(j)$ at $j$. Note that $U_B(\mu_1+\mu_2)=U_B(\mu_1)+U_B(\mu_2)$ if the supports of $\mu_1$ and $\mu_2$ are disjoint. For any trajectory $\omega$ denote by $U_{B,t}(\omega)=U_B(\omega_t)$ that object at time $t$.

Since the distribution of $Y_t$, $\Poy_{t,\rho}$, is infinitely divisible, $U_{B,t}$ is a Poisson process.

\lemma{ \label{thm:mprw:poisson}
For $t\in[0,1)$, let $U_{B,t}$ be the measure counting the number of points according to their multiplicity in a given set $B\in\Bbd$. Then $U_{B,t}$ is a finite Poisson process with intensity measure
\begin{equation*}
	\rho(B)\sum_{j\geq 1} \frac{t^j}{j}\delta_j.
\end{equation*}
}
In particular the total mass of $U_{B,t}$, $U_{B,t}(\N)$, is the total number of occupied sites of $Y_t$ in $B$, and is Poisson distributed with intensity $-\rho(B)\log(1-t)$. Denote by $\tau_t$ the measure $\sum_{j} \frac{t^j}{j}\delta_j$.

What about the dynamics of the process $(U_{B,t})_t$ in finite sets $B$? Note that by the convolution property of the P\'olya sum process the transition kernel satisfies for $s<t$
\begin{align}
	p_{s,t}(\nu,\phi)
		&=\int \phi(\mu+\nu) \Poy_{\frac{t-s}{1-s},\rho+\nu}(\d\mu) \notag\\
		&=\iint \phi(\mu_i+\mu_b+\nu)\Poy_{\frac{t-s}{1-s},\rho}(\d\mu_i) \Poy_{\frac{t-s}{1-s},\nu}(\d\mu_b).\label{eq:mprw:scrp-trans}
\end{align}
Thus given a point configuration or population $\nu$, the transition from $\nu$ to a new configuration means to add two parts, in terms of a population an 'immigration' configuration and an 'offspring' configuration.

More precisely, the transition can be described as follows: During the time intervall $(s,t]$, new points may immigrate, they are allowed to arrive as multiple points. Locational preference depends on $\rho$. Since we assumed $\rho$ to be diffuse, they do not hit already occupied locations. Independently, already existing points may branch and share their location with their offsprings.

\lemma{
For a bounded set $B\in\B$ and diffuse measure $\rho$, $(U_{B,t})_t$ is a Markov process on the state space $\Mpm(\N)$ with transition kernel
\begin{equation*}
	p^B_{s,t}(\eta,f)=\int f(\eta_i+\eta_b)\bigl(U_B\Poy_{\frac{t-s}{1-s},\rho}\bigr)(\d\eta_i) \bigl(U_B(\Poy_{\frac{t-s}{1-s},\nu}\ast\Delta_\nu)\bigr)(\d\eta_b)
\end{equation*}
for a.e. $\nu\in\{U_B=\eta\}$.
}

\begin{proof}
First note that $U_B:\MpmX\to\MpfN$ is onto. Secondly, set $\phi=f\circ U_B$ in~\eqref{eq:mprw:scrp-trans} for some measurable and bounded function $f:\Mpf(\N)\to\R$, then
\begin{align*}
	p_{s,t}(\nu,f\circ U_B) 
		&= \iint f\bigl(U_B(\mu_i+\mu_b+\nu)\bigr) \Poy_{\frac{t-s}{1-s},\rho}(\d\mu_i) \Poy_{\frac{t-s}{1-s},\nu}(\d\mu_b)\\
		&= \iint f\bigl(U_B(\mu_i)+U_B(\mu_b+\nu)\bigr) \Poy_{\frac{t-s}{1-s},\rho}(\d\mu_i) \Poy_{\frac{t-s}{1-s},\nu}(\d\mu_b)
\end{align*}
since $\rho$ is diffuse and therefore $\supp\mu_i$ and $\supp(\mu_i+\nu)=\supp\nu$ are disjoint. Note that if $\nu'$ is such that $U_B\nu=U_B\nu'$, then $U_B(\mu_b+\nu)$ and $U_B(\mu_b+\nu')$ have the same distribution under $\Poy_{\frac{t-s}{1-s},\nu}$ and $\Poy_{\frac{t-s}{1-s},\nu'}$, respectively, and $p_{s,t}(\ccdot,f\circ U_B)$ is a.s. constant on sets $\{U_B=\eta\}$.
\end{proof}


\korollar{ \label{thm:mprw:gen}
For a bounded set $B\in\B$, the generator of $(U_{B,t})_t$ is given by
\begin{equation*}
	\A^B_s f(\eta) = \frac{\rho(B)}{1-s}\bigl[f(\eta+\delta_1)-f(\eta)\bigr] +\frac{1}{1-s} \sum_{j\geq 1} j\eta(j)\bigl[f(\eta-\delta_j+\delta_{j+1})-f(\eta)\bigr]
\end{equation*}
for $f:\MpfN\to\R$ measurable and bounded.
}
The series converges since the point measure $\eta$ has finite mass and therefore finite first moment.

\begin{proof}
Let $\eta\in\MpfN$ and $\mu\in\{U_B=\eta\}$, then
\begin{align*}
	\A U_B^\ast f(\mu)
		&=\frac{1}{1-s}\int_B f\bigl(U_B(\mu+\delta_x)\bigr)-f(U_B\mu)\bigl(\rho+\mu\bigr)(\d x)\\
		&=\begin{multlined}[t]
				\frac{1}{1-s}\int_B f(\eta+\delta_1)-f(\eta)\rho(\d x) \\
					+\frac{1}{1-s}\int_B f\bigl(U_B(\mu+\delta_x)\bigr)-f(U_B\mu)\mu(\d x)
			\end{multlined}\\
		&=\frac{\rho(B)}{1-s} \bigl[ f(\eta+\delta_1)-f(\eta) \bigr] 
			+\frac{1}{1-s}\sum_{j\geq 1} j\eta(j) \bigl[f(\eta-\delta_j+\delta_{j+1})-f(\eta)\bigr]
\end{align*}
since a stack of $j$ individuals may give each of its particles the chance to branch.
\end{proof}

\bem{
$\A^B_s$ may computed directly along similar lines as $\A_s$.
}

This representation allows a different point of view on the growth process $(U_{B,t})_t$: Observing a single fixed site, as time passes, particles arrive from the left and leave to the right. Meanwhile the stack grows in mean. From the point of view of some moving particle, from times $s$ to time $t$  particle at site $j$ steps to the right with a $NB_{\frac{t-s}{1-t},j}$-distributed step independently of other particles. Moreover, between times $s$ and $t$, new sites get occupied. Firstly remark that this number is Poisson, and secondly that each of the new particles on $\N$ gets a random starting point assigned, which is chosen geometrically distributed independent of all other mechanisms. The surprising fct is that nevertheless at some fxed time, the counting variables of the number of particles at some site are independent.

Finally we compute the generator of the centered process $(U_{B,t}-\rho(B)\tau_t)_t$, that is for each time $t$ a finite signed random measure on $\N$. Since this is not a state space transformation covered by the state space transformation theorem, we compute the generator directly. A direct proof of the claim in Corollary~\ref{thm:mprw:gen} could hve been obtained in a very similar but simpler fashion.

For a measurable function $f:\MpfN\to\R$ denote by $f':\MffN\times\MffN\to\R$, $(\eta,\xi)\mapsto f'(\eta)[\xi]$, its Gateaux derivative at $\eta$ in direction $\xi$ in case of existence.

\prop{ \label{thm:mprw:gen-centered}
For a bounded set $B\in\B$, the generator of $(U_{B,t}-\rho(B)\tau_t)_t$ is given by
\begin{equation*}
	\C^B_s f(\xi) = \begin{multlined}[t]
			\frac{\rho(B)}{1-s}\bigl[f(\xi+\delta_1)-f(\xi)\bigr] -\rho(B) \sum_{j\geq 1} s^{j-1}f'(\xi)[\delta_j] \\
			+\frac{1}{1-s} \sum_{j\geq 1} \bigl(j\xi(j)+\rho(B)s^j\bigr)\bigl[f(\eta-\delta_j+\delta_{j+1})-f(\eta)\bigr]
		\end{multlined}
\end{equation*}
for $f:\MffN\to\R$ continuously differentiable with bounded derivative and $\xi\in\MffN$.
}

\begin{proof}
Mainly the arguments are analogue to those of the proof of Lemma~\ref{thm:gen:gen}. First note that the transition kernels are given by
\begin{equation*}
	q^B_{s,t}(\xi,f)=p^B_{s,t}\bigl(\xi+\rho(B)\tau_s,f(\ccdot -\rho(B)\tau_t)\bigr).
\end{equation*}
Then for $\eta=\xi+\rho(B)\tau_s$ and some $\nu\in\{U_B=\eta\}$,
\begin{equation*}
	\begin{multlined}[t]
			q^B_{s,s+h}(\xi,f)-f(\xi)\\ 
	= \left[\left(1-\frac{h}{1-s}\right)^{\rho(B)+\nu(B)}-1\right]f\bigl(\eta-\rho(B)\tau_{s+h}\bigr)\\
				+f\bigl(\eta-\rho(B)\tau_{s+h}\bigr)-f\bigl(\eta-\rho(B)\tau_{s}\bigr)\\
				+ \int_{\zeta_B\geq 1}\int_{\zeta_B=0} f\bigl(U_B(\mu_i)+U_B(\mu_b+\nu)-\rho(B)\tau_{s+h}\bigr) \Poy_{\frac{h}{1-s},\rho}(\d\mu_i)\Poy_{\frac{h}{1-s},\nu}(\d\mu_b)\\
				+ \int_{\zeta_B=0}\int_{\zeta_B\geq 1} f\bigl(U_B(\mu_i)+U_B(\mu_b+\nu)-\rho(B)\tau_{s+h}\bigr) \Poy_{\frac{h}{1-s},\rho}(\d\mu_i)\Poy_{\frac{h}{1-s},\nu}(\d\mu_b)\\
				+ \int_{\zeta_B\geq 1}\int_{\zeta_B\geq 1} f\bigl(U_B(\mu_i)+U_B(\mu_b+\nu)-\rho(B)\tau_{s+h}\bigr) \Poy_{\frac{h}{1-s},\rho}(\d\mu_i)\Poy_{\frac{h}{1-s},\nu}(\d\mu_b),
			\end{multlined}
\end{equation*}
where $\nu\in\{U_B=\eta\}$. Divide by $h$ and let $h\searrow 0$, then the first summand tends to $-\tfrac{\rho(B)+\nu(B)}{1-s}f(\xi)=-\tfrac{\rho(B)+J\eta}{1-s}f(\xi)$, where $J\eta=\sum_j j\eta(j)=J\xi+\sum_j s^j\delta_j$. For the last two summands, the probability of realizing no point tends to 1, and the remaining immigration part to $\tfrac{1}{1-s}\int f\bigl(U_B(\nu+\delta_x)-\rho(B)\tau_{s+h}\bigr)\rho(\d x)= \tfrac{\rho(B)}{1-s} f(\xi+\delta_1)$ since $\rho$ is diffuse, but the branching part to $\tfrac{1}{1-s}\int f\bigl(U_B(\nu+\delta_x)-\rho(B)\tau_{s+h}\bigr)\nu(\d x)$. Note that this integral does only depend on $\nu$ via $U_B(\nu)=\eta$, more precisely,
\begin{equation*}
	\frac{1}{1-s}\int f\bigl(U_B(\nu+\delta_x)-\rho(B)\tau_{s+h}\bigr)\nu(\d x) = \frac{1}{1-s} \sum_{j\geq 1} j\eta(j)f(\xi-\delta_j+\delta_{j+1}).
\end{equation*}
Finally, the last summand vanishes. 
\end{proof}

\section{Hydrodynamic limit \label{sect:hydra}}

If $\rho$ is an infinite measure, as the observation window $B$ grows, the number ob observed particles grows unboundedly. Thus the first step is to exhaust $X$ and renormalize that $U_B$.

Let $V_{B,t}=\tfrac{U_{B,t}}{\rho(B)}$ be the renormalization, then $(V_{B,t})_t$ is a Markov process.
\lemma{
$(V_{B,t})_{0\leq t<1}$ is a Markov process with transition kernel
\begin{equation*}
	\tilde{p}^B_{s,t}(\eta,f)
			=\int f\left(\frac{\xi}{\rho(B)}\right) p^B_{s,t}\bigl(\rho(B)\eta,\d\xi\bigr)
\end{equation*}
and generator
\begin{equation} \label{eq:hydro:gen:norm}
	\tilde{\A}^B_s f(\eta) = \begin{multlined}[t]
				\frac{\rho(B)}{1-s}\left[f\left(\eta+\frac{\delta_1}{\rho(B)}\right)-f(\eta)\right]\\
			 +\frac{\rho(B)}{1-s} \sum_{j\geq 1} j\eta(j)\left[f\left(\eta-\frac{\delta_j+\delta_{j+1}}{\rho(B)}\right)-f(\eta)\right].
			 \end{multlined}
\end{equation}
}
Let $(B_n)_n$ an increasing sequence of bounded and measurable sets which exausts $X$. Then the law of lage numbers yields
\prop{
Let $\rho$ be an infinite and diffuse measure, then for each $t\in[0,1)$ almost surely
\begin{equation*}
	V_t=\lim_{n\to\infty} V_{B_n,t} = \sum_{j\geq 1} \frac{t^j}{j}\delta_j.
\end{equation*}
}
We denote by $\tau_t$ the limitng measure on the right hand side.
\prop{
This convergence holds for the process $(V_{B_n,t})_t$
}
\begin{proof}
Choosing a countable, dense subset of $[0,1)$ to get a joint exceptional set, the right continuity of the paths yields the convergence.
\end{proof}

Thus the process converges, and what remains is the identification of the dynamics, which should describe a flow of mass on the positive integers from 1 to the right. 

\korollar{
Let $f:\MffN\to\R$ be twice continuously Gateaux differentiable with bounded second derivative and $\tilde{\A}$ be the generator of $(V_t)_{0\leq t<1}$. Then $f\in\dom\tilde{\A}$ and
\begin{equation*}
	\tilde{\A}_s f(\eta) =\frac{1}{1-s} \left[ f'(\eta)[\delta_1]+\sum_{j\geq 1} j\eta(j) f'(\eta)[-\delta_j+\delta_{j+1}]\right].
\end{equation*}
}
\begin{proof}
Let $f$ be twice continuously Gateaux differentiable with bounded second derivative, then
\begin{equation*}
	f(\eta+r\kappa)-f(\eta)=f'(\eta)[r\kappa]+o(r)=r f'(\eta)[\kappa]+ o(r).
\end{equation*}
Thus inserting this into~\eqref{eq:hydro:gen:norm} and letting $n\to\infty$ yields the result.
\end{proof}

\bem{
\begin{enumerate}
	\item Observe that explicitly
\begin{align*}
	\tilde{\A}_sf(\tau_s) &= \frac{1}{1-s}\left[ f'(\tau_s)[\delta_1]+\sum_{j\geq 1}s^j f'(\tau_s)[-\delta_j+\delta_{j+1}] \right]\\
		&=\sum_{j\geq 1} s^{j-1}f'(\tau_s)[\delta_j] =	\frac{\d\phantom{s}}{\d s}f(\tau_s),
\end{align*}
which is explicitly Kolmogorow's equation for the hydrodynamic limit $(V_t)$.
	\item Rearrangement of the terms in the generator yiels
		\begin{equation*}
			\tilde{\A}_s f(\eta) =\frac{1}{1-s} \left[ \bigl(1-\eta(1)\bigr)f'(\eta)[\delta_1]+\sum_{j\geq 2} \bigl[(j-1)\eta(j-1)-j\eta(j)\bigr] f'(\eta)[\delta_j] \right]
		\end{equation*}
		and we get explicitly the rates at which the mass being introduced at 1 moves to the right. More explicitly, by testing successively $f=\zeta_{\{j\}}$, for which $\zeta_{\{j\}}=1\otimes\zeta_{\{j\}}$, we get
		\begin{align*}
			\frac{\d\phantom{t}}{\d t} V_t(j) = \begin{cases}
					\frac{1}{1-t}\bigl[1-V_t(1)\bigr]  & j=1,\\
					\frac{1}{1-t}\bigl[(j-1)V_t(j-1)-jV_t(j)\bigr]  & j>1.
				\end{cases}
		\end{align*}
\end{enumerate}
}

\section{Fluctuations}

In the similar manner we want to identify the fluctuations of the process and how they behave. By Lemma~\ref{thm:mprw:poisson} we know that for each $j$, $U_{B,t}(j)$ is Poisson distributed with mean $\tfrac{t^j}{j}\rho(B)$, and that $U_{B,t}(1),U_{B,t}(2),\ldots$ are independent random variables. Let
\begin{equation*}
	Z_{B,t}:\MpfN\to\MN,\qquad 
	Z_{B,t}(\eta)=\frac{1}{\sqrt{\rho(B)}}\sum_{j\geq 1} \left[\eta(j)-\frac{t^j\rho(B)}{j}\right]. 
\end{equation*}
with abuse of notation denote in the following by $Z_{B,t}$ also the random variable $Z_{B,t}\circ U_{B,t}$. Then, as $B$ increases, the central limit theorem ensures normal behaviour.
\lemma{ \label{thm:fluc:clt}
Let $(B_n)_n$ be an increasing sequence of bounded sets exhausting $X$. Then for each $t\in(0,1)$ and $j\in\N$, $Z_{B_n,t}(j)$ converges weakly to a $\Nd(0,\tfrac{t^j}{j})$-distributed random variable $Z_t(j)$. Moreover, the random variables $Z_t(1),Z_t(2),\ldots$ are independent.
}
 
A direct consequence of Dynkin's state space transformation theorem in connection with Proposition~\ref{thm:mprw:gen-centered} shows that $(Z_{B,t})_t$ is a Markov process.

\prop{
$(Z_{B,t})_{0\leq t< 1}$ is a Markov process, and if $(\C^B_s)_{0\leq s<1}$ denotes the generating family of its tansition semigroup, then for $f:\MffN\to\R$ continuously differentiable with bounded derivative and $\eta\in\MffN$,
\begin{equation*}
	\begin{multlined}[t]
		\C^B_{s}f(\eta) =\frac{\rho(B)}{1-s} \left[f\left(\eta+\frac{\delta_1}{\sqrt{\rho(B)}}\right)-f(\eta)\right]-\sqrt{\rho(B)} \sum_{j\geq 1} s^{j-1}f'(\xi)[\delta_j]\\
		+\frac{1}{1-s}\sum_j \left(\sqrt{\rho(B)}j\xi(j)+\rho(B)s^j\right) \left[f\left(\eta-\frac{\delta_j-\delta_{j+1}}{\sqrt{\rho(B)}}\right)-f(\eta)\right].
	\end{multlined}
\end{equation*}
}

Again, let $(B_n)_n$ be an increasing sequence of bounded sets which exhausts $X$ as well as $f$ a three times Gateaux differentiable function with bounded third derivative. Then
\begin{align*}
	&f(\xi+\sqrt{r}\delta_1)-f(\xi) 
		= \sqrt{r}f'(\xi)[\delta_1]+\frac{r}{2}f''(\xi)[\delta_1,\delta_1] +o(r)\\
	&\begin{multlined}[t] 
		f(\xi-\sqrt{r}\delta_j+\sqrt{r}\delta_{j+1})-f(\xi) \\
		= \sqrt{r} \left[ f'(\xi)[\delta_{j+1}]-f'(\xi)[\delta_j] \right]\\
			+ \frac{r}{2}\left[ f''(\xi)[\delta_{j+1},\delta_{j+1}]+f''(\xi)[\delta_j,\delta_j] -f''(\xi)[\delta_j,\delta_{j+1}]\right] +o(r)
	\end{multlined}
\end{align*}
Inserting these with $r$ replaced by $\rho(B_n)^{-1}$ into $\C^B_s$ and letting $n\to\infty$ yields
\prop{
Let $Z_t=\lim_{n\to\infty} Z_{B_n,t}$ and denote by $(\C_s)_{0\leq s<1}$ the generating fmily of $(Z_t)_{0\leq t<1}$. Then if $f:\MffN\to\R$ is three times differentiable with bounded derivative and $\xi\in\MffN$,
\begin{equation*}
	\C_s f(\xi)=\begin{multlined}[t]
		\frac{1}{2}\cdot\frac{1+s}{1-s}\sum_{j\geq 1} s^{j-1} f''(\xi)[\delta_j,\delta_j]
			-\frac{1}{1-s}\sum_{j\geq 1} s^j f''(\xi)[\delta_j,\delta_{j+1}]\\
			+\frac{1}{1-s}\sum_{j\geq 1} \left[f'(\xi)[\delta_{j+1}]-f'(\xi)[\delta_j]\right]j\xi(j)
		\end{multlined}
\end{equation*}
}
\begin{proof}
What remains is to expand, to sort and to let $n\to\infty$. Observe that the $\sqrt{\rho(B_n)}$-terms are
\begin{equation*}
	\begin{multlined}[t]
		\frac{\sqrt{\rho(B_n)}}{1-s}f'(\xi)[\delta_1]+\frac{\sqrt{\rho(B_n)}}{1-s}\sum_{j\geq 1} s^j\left[f'(\xi)[\delta_{j+1}]-f'(\xi)[\delta_j]\right]\\
		=\sqrt{\rho(B_n)}\sum_{j\geq 1}s^{j-1}f'(\xi)[\delta_j],
	\end{multlined}
\end{equation*}
which cancels the derivative term in $\C^B_sf$. Expanding the second derivative terms and collecting the pure derivatives gives a $\tfrac{1+1}{1-s}$
\end{proof}

\bem{
A rearrangement of the drift term yields
\begin{equation*}
	\C_s f(\xi)=\begin{multlined}[t]
		\frac{1}{2}\cdot\frac{1+s}{1-s}\sum_{j\geq 1} s^{j-1}f''(\xi)[\delta_j,\delta_j]
			-\frac{1}{1-s} \sum_{j\geq 1} s^j f''(\xi)[\delta_j,\delta_{j+1}]\\
			-\frac{1}{1-s}\sum_{j\geq 1} \left[j\xi(j)-(j-1)\xi(j-1)\right] f'(\xi)[\delta_j],
		\end{multlined}
\end{equation*}
where we put $\xi(0)=0$ for convenience. Thus at each $j$ there is a Brownian motion with some drift being influenced by the left neighbouring Brownian motion and some negative cross diffusion coefficient.
}

With the aid of moment generating functions, i.e. for functions $f$ of type $\e^{\zeta_g}$ with a bounded function $g:\N\to\R$, we finally compute $\Ex\C_sf(Z_s)$. Note that $f'=f\otimes g$ and $f''=f\otimes g\otimes g$. Moreover, since for $Z\sim\Nd(0,\sigma^2)$,
\begin{align*}
	\Ex\bigl[Z\e^{uZ}\bigr] &= \frac{1}{u}\frac{\d\phantom{t}}{\d u}\Ex\bigl[\e^{uZ}\bigr]\Big|_{t=1}
		=u\sigma^2\Ex\bigl[\e^{uZ}\bigr],
\end{align*}
we get $\Ex\bigl[Z_s(j)\e^{Z_s(g)}\bigr]=g(j)\tfrac{s^j}{j}\Ex\e^{Z_s(g)}$. Thus
\begin{align*}
	\Ex\C_sf(Z_s) &= \Ex\bigl[\e^{uZ}\bigr]
		\begin{multlined}[t]
			\left\{ \frac{1}{2}\cdot\frac{1+s}{1-s}\sum_{j\geq 1} s^{j-1} g(j)^2-\frac{1}{1-s}\sum_{j\geq 1}s^jg(j)g(j+1)\right.\\ 
			\left.+\frac{1}{1-s}\sum_{j\geq 1}s^jg(j)g(j+1)-\frac{1}{1-s}\sum_{j\geq 1}s^jg(j)^2 \right\}
		\end{multlined}\\
		&=\Ex\bigl[\e^{uZ}\bigr]\cdot\frac{1}{2}\sum_{j\geq 1} s^{j-1} g(j)^2,
\end{align*}
which agrees with
\begin{equation*}
	\frac{\d\phantom{s}}{\d s}\Ex\bigl[Z\e^{uZ}\bigr] 
		= \frac{\d\phantom{s}}{\d s}\exp\left[\sum_{j\geq 1}\frac{s^j}{j}g(j)^2\right]
\end{equation*}
according to Lemma~\ref{thm:fluc:clt}.

\end{document}